\documentclass[11pt]{amsart}
\usepackage{amsmath,amsthm, amscd, amssymb, amsfonts}

\usepackage{epsfig}
\usepackage[all]{xy}

\newcommand\boxe{\begin{tabular}{|p{0,1cm}|}
\hline \\ \hline \end{tabular}}

\newcommand{\sou}{\mathfrak s}
\newcommand{\tgt}{\mathfrak e}
\DeclareMathOperator*{\Tim}{\times}
\newcommand{\ftimes}{\sideset{_\tgt}{_\sou}\Tim}
\newcommand{\hvtimes}{\sideset{_r}{_t}\Tim}
\newcommand{\dtimes}{\sideset{_b}{_l}\Tim}

\newcommand{\Vc}{{\mathcal V}}
\newcommand{\Hc}{{\mathcal H}}
\newcommand{\B}{{\mathcal B}}
\newcommand{\Bg}{{\mathfrak B}}

\newcommand{\clase}{{\mathfrak{U}}}
\newcommand{\clate}{{\mathfrak{R}}}

\newcommand{\Mo}{{\mathcal M}}
\newcommand{\G}{{\mathcal G}}
\newcommand{\D}{{\mathcal D}}
\newcommand{\Ec}{{\mathcal E}}
\newcommand{\Pc}{{\mathcal P}}
\newcommand{\X}{{\mathcal X}}

\newcommand{\fde}{{\triangleright}}
\newcommand{\gde}{{\triangleleft}}
\newcommand{\kuG}{\ku \, {\mathcal G}}
\newcommand{\kuH}{\ku \, {\mathcal H}}
\newcommand{\Tc}{ {\mathcal T}}
\newcommand{\moh}{\,_H{\mathcal M}}
\newcommand{\mog}{\,_{\kk\G}{\mathcal M}}
\newcommand{\moha}{\,_{\kuH}{\mathcal M}}
\newcommand{\mogp}{\,_{\kuG (P)}{\mathcal M}}

\newcommand{\Ss}{{\mathcal S}}

\newcommand{\Ker}{\mbox{\rm Ker\,}}

\newcommand{\Do}{\Vc\bowtie\Hc}

\newcommand{\Qc}{{\mathcal Q}}
\newcommand{\s}{{\sigma}}

\newcommand{\kk}{{\mathbb K}}
\newcommand{\ku}{{\Bbbk}}
\newcommand{\Z}{{\mathbb Z}}
\newcommand{\Na}{{\mathbb N}}

\newcommand{\Sy}{{\mathbb S}}

\newcommand{\C}{{\mathbb C}}

\newcommand\id{\operatorname{id}}
\newcommand{\Tot}{\mbox{\rm Tot\,}}

\newcommand{\Res}{\mbox{\rm Res\,}}

\newcommand{\Aut}{\mbox{\rm Aut\,}}

\newcommand{\Hom}{\mbox{\rm Hom\,}}

\newcommand{\Ext}{\mbox{\rm Ext\,}}
\newcommand{\coH}{\mbox{\rm H\,}}

\newcommand\Opext{\operatorname{Opext}}

\newcommand{\im}{\mathop{\rm Im\,}}

\newcommand{\Id}{\mathop{\rm Id}}

\renewcommand{\Im}{\im}

\theoremstyle{plain}

\numberwithin{equation}{section}

\newtheorem{teo}{Theorem}[section]

\newtheorem{lema}[teo]{Lemma}

\newtheorem{cor}[teo]{Corollary}

\newtheorem{prop}[teo]{Proposition}

\theoremstyle{definition}

\newtheorem{defi}[teo]{Definition}
\newtheorem{exa}[teo]{Example}
\newtheorem{sm}[teo]{Summary}

\theoremstyle{remark}

\newtheorem{rmk}[teo]{Remark}

\def\pf{\begin{proof}}

\def\epf{\end{proof}}

\theoremstyle{remark}

\begin{document}

\title[Examples of face algebras]{Examples of weak Hopf algebras
arising from vacant double groupoids}
\author[Andruskiewitsch and Mombelli]{Nicol\'as Andruskiewitsch and  Juan
Mart\'\i n Mombelli}
\thanks{This work was partially supported by
Agencia C\'ordoba Ciencia, ANPCyT-Foncyt, CONICET, ECOS,
Fundaci\'on Antorchas and Secyt (UNC)}
\address{Facultad de Matem\'atica, Astronom\'\i a y F\'\i sica
\newline \indent
Universidad Nacional de C\'ordoba
\newline
\indent CIEM -- CONICET
\newline \indent Medina Allende s/n
\newline
\indent (5000) Ciudad Universitaria, C\'ordoba, Argentina}
\email{andrus@mate.uncor.edu, \quad \emph{URL:}\/
http://www.mate.uncor.edu/andrus} \email{mombelli@mate.uncor.edu}
\begin{abstract} We construct explicit examples of weak Hopf algebras (actually
face algebras in the sense of Hayashi \cite{H}) via vacant double
groupoids as explained in \cite{AN}. To this end, we first study
the Kac exact sequence for matched pairs of groupoids and show
that it can be computed via group cohomology. Then we describe
explicit examples of finite vacant double groupoids.
\end{abstract}

\date{\today}
\maketitle

\section*{Introduction}

Tensor categories have many important applications in several
areas of mathematics and theoretical physics. A source of examples
of tensor categories is the theory of Hopf algebras; namely the
category of representations of a Hopf algebra is naturally a
tensor category. However, there are important linear tensor
categories that do not arise as the category of representations of
any Hopf algebra. Some fifteen years ago, Ocneanu proposed the
notion of paragroup to encompass these examples. In this
direction, Hayashi introduced face algebras in 1991; eventually,
he showed that a suitable linear tensor category arises as the
category of representations of a face algebra (canonical Tannaka
reconstruction). See \cite{H} and references therein. Weak Hopf
algebras were introduced in \cite{bnsz, bsz}; face algebras are
weak Hopf algebras with commutative target subalgebra.

Recently, it was explained how to build up weak Hopf algebras
(actually face algebras) starting from a matched pairs of finite
groupoids and a suitable pair of cocycles \cite{AN}. The purpose
of the present paper is to exhibit explicit examples of weak Hopf
algebras in the framework of this construction. For this, we need
to give explicit examples of matched pairs of finite groupoids
(what we do in Section 3), and to compute the corresponding 2
cohomology group (the so-called Opext). As said in \cite{AN}, see
also \cite{BSV}, an efficient way for this last task is through
the Kac exact sequence, a generalization of the analogous sequence
for matched pairs of groups. We elaborate on Kac exact sequences
in Section 2, relating to cohomology of weak Hopf algebras, which
we discuss in Section 1.

The reader interested in the construction of explicit examples
might find useful the Summaries \ref{genconst}, \ref{case1},
\ref{caso2} and \ref{case3}. We include along the way some
calculations of the Opext groups, by reduction to group
cohomology.

\subsection*{Notation} We shall denote by $\kk$ a commutative ring and
by $\ku$ a field of characteristic zero. If $R$ is an algebra, we
denote by $_R\Mo$ the category of left $R$-modules. If $\X$ is a
set, we denote by $\kk \, \X$ the free $\kk$-module with basis
$(X)_{X\in \X}$. We shall use Sweedler's notation but omitting the
summation sign for coalgebras: $\Delta(x) = x_{(1)} \otimes
x_{(2)}$, if $\Delta$ is the comultiplication of a coalgebra $C$,
$x\in C$.
For any ring $R$ we shall denote by $R^{\times}$ the group of
invertible elements in $R$.

\subsection*{Acknowledgments} We thank Sonia Natale for many
conversations. Part of the work of the second author was
done during a visit to the University of Rheims in the framework
of the project ECOS. He is very grateful to Jacques Alev for his
kind hospitality.

\bigbreak
\section{Cohomology of groupoids}

\subsection{Weak Hopf Algebras}

\

We first recall the definition of weak Hopf algebras, or quantum
groupoids \cite{bnsz, bsz}; see also \cite{NV}.
A \emph{weak bialgebra} is a collection $(H, m, \Delta)$, where
$(H,m)$ is an associative $\kk$-algebra with unit 1 and $(H,
\Delta)$ is a coassociative $\kk$-coalgebra  with counit
$\varepsilon$, such that the following axioms hold:
\begin{align}\label{d-mult}\Delta(ab) &= \Delta(a) \Delta(b), \qquad
\forall a, b \in H.
\\ \label{ax-unit} \Delta^{(2)} (1 ) &= \left( \Delta(1) \otimes 1 \right)
\left( 1 \otimes \Delta(1) \right) = \left( 1 \otimes \Delta(1)
\right) \left( \Delta(1) \otimes 1 \right).
\\ \label{ax-counit} \varepsilon(abc) &= \varepsilon(ab_{(1)})
\varepsilon(b_{(2)}c) =
\varepsilon(ab_{(2)})\varepsilon(b_{(1)}c), \qquad \forall a, b, c
\in H.
\end{align}

A weak bialgebra $H$ is a \emph{weak Hopf algebra} or a
\emph{quantum groupoid} if there exists a linear map $\Ss: H \to
H$ satisfying
\begin{align}\label{atp-1} m(\id \otimes \Ss) \Delta (h) &
= (\varepsilon \otimes \id) \left(  \Delta(1) (h \otimes 1)\right)
=:
\varepsilon_t(h), \\
\label{atp-2} m(\Ss \otimes \id) \Delta (h) & = (\id \otimes
\varepsilon) \left( (1 \otimes h) \Delta(1) \right)=:
\varepsilon_s(h),\\
\label{atp-3} m^{(2)}(\Ss \otimes \id \otimes \Ss) \Delta^{(2)} &
= \Ss,
\end{align}
for all $h \in H$. The maps $\varepsilon_s$, $\varepsilon_t$ are
respectively called the source and target maps; their images are
called the source and target subalgebras, and we denote them
respectively by $H_s$ and $H_t$. 

The weak Hopf algebra $H$ is an $H_t$-bimodule via $z.h.w:=zhw$
for $h\in H$, $z,w\in H_t$ and the target subalgebra $H_t$ has a
left $H$-module structure given by:
$$h.w=\epsilon_t(hw),$$
for all $h\in H,$ and $w\in H_t$. This action when restricted to
$H_t$ gives the left regular action \cite[p. 215]{NV}.
The following Lemma will be useful later.
\begin{lema}\label{tec} Let be $H$ a weak Hopf algebra, $H_t$ its target
subalgebra, $M$ a left $H_t$-module and $N$ a left $H$-module.
Then $H\otimes_{H_t} M$ has a left $H$-module structure via
multiplication on the first tensorand and there are natural
isomorphisms
$$\Hom_H(H\otimes_{H_t} M,N)\;\simeq\; \Hom_{H_t}(M,\Res_{H_t}^{H} N).$$

For any projective $H_t$-module $M$, the $H$-module
$H\otimes_{H_t} M$ is projective.

\end{lema}

\begin{proof} The desired natural isomorphisms are defined by
\begin{align*} \phi:\Hom_H(H\otimes_{H_t} M,N)\rightarrow
\Hom_{H_t}(M,\Res_{H_t}^{H} N),\;\;
\phi(f)(m)=f(1\otimes m), \\
\psi:\Hom_{H_t}(M,\Res_{H_t}^{H} N)\rightarrow
\Hom_H(H\otimes_{H_t} M,N),\;\; \psi(g)(h\otimes m)=h g(m),
\end{align*}
for all $h\in H$, $m\in M$. The last claim follows from the first
one.
\end{proof}

\begin{rmk} If $\kk = \ku$ is a field then
by \cite[Prop. 2.3.4]{NV} the target subalgebra is separable and
therefore semisimple and thus every $H_t$-module is projective. In
this case, if $H_t$ is commutative then $H$ is a face algebra in
the sense of Hayashi \cite{H}.

In fact, it seems that weak Hopf algebras are always considered
over a field, and we were not able to find a thorough study of
weak Hopf algebras over more general commutative rings in the
literature. However, groupoid algebras over $\kk$ are examples of
such objects.
\end{rmk}

\medbreak
\subsection{ The Bar Resolution for Weak Hopf Algebras}\label{weak}

\

Let $H$ be a weak Hopf algebra with target subalgebra $H_t$. We
define the cohomology groups of $H$ with coefficients in $M\in
\moh$ by
$$\coH^n(H,M):=\Ext^n_H(H_t,M).$$

These cohomology groups can be computed by means of a ``normalized
bar resolution". Let $\overline{H}$ be the $H_t$-bimodule $H/H_t$;
and let $\overline{h}$ be the class in $\overline{H}$ of $h\in H$.
If $N \in \moh$ and $n\in {\Na}_0$, we set
$$B_n(H,N):=H\otimes_{H_t}\underbrace{\overline{H}\otimes_{H_t}...\otimes_{H_t}
\overline{H}}_{n-times}\otimes_{H_t}N.$$

Then $B_n(H,N)$ is a left $H$-module via multiplication on the
first tensorand, and if $N$ is a projective $H_t$-module then
$B_n(H,N)$ is a projective $H$-module thanks to Lemma \ref{tec}.
Let us define maps $\epsilon:B_0(H,N)\rightarrow N$,
$\partial_n:B_n(H,N)\rightarrow B_{n-1}(H,N)$, $n > 0$, and
$s_n:B_n(H,N)\rightarrow B_{n+1}(H,N)$, $n\geq 0$,  by
\begin{align*}
\epsilon(h\otimes m)&=h.m, \\
\partial_n(h\otimes \overline{h_1}\otimes
...\otimes\overline{h_n}\otimes m)&=hh_1\otimes \overline{h_2}
\otimes...\otimes \overline{h_n}\otimes m\\ + \sum_{i=1}^{n-1}
(-1)^i h\otimes \overline{h_1}\otimes ... \otimes \overline{h_i
h_{i+1}}\otimes...\otimes \overline{h_n}\otimes m  &+ (-1)^n
h\otimes \overline{h_1}\otimes...\otimes\overline{h_{n-1}}\otimes
h_n. m, \;\;
\end{align*}
and
$$s_n(h\otimes
\overline{h_1}\otimes...\otimes\overline{h_n}\otimes m)=1\otimes
\overline{h}\otimes\overline{h_1}\otimes...\otimes\overline{h_n}\otimes
m,$$

\noindent for all $h, h_1,\dots ,h_n\in H$ and for all $m\in N$.

\begin{lema}\label{tec3}For any $n > 0$, we have
\begin{itemize}
\item[i)]  the maps $\partial_n$ are well defined $H$-module
homomorphisms, \item[ii)] $\partial_{n+1}\partial_n=0$, and
\item[iii)] $\partial_{n+1} s_n + s_{n-1}
\partial_n=\id_{B_n(H,N)}.$
\end{itemize}
\end{lema}
\begin{proof}  We verify i). Assume that $h_i\in H_t$ for $1<i<n$ then
\begin{align*}\partial_n(h\otimes
\overline{h_1}\otimes...&\otimes\overline{h_n} \otimes
m)=(-1)^{i-1} h\otimes ... \otimes \overline{h_{i-1}
h_{i}}\otimes...\otimes m +(-1)^i h\otimes ... \otimes
\overline{h_i
h_{i+1}}\otimes...\otimes m\\
&=(-1)^{i-1} h\otimes ... \otimes \overline{h_{i-1} }\otimes
h_{i}\overline{h_{i+1}} ...\otimes m +(-1)^i h\otimes ... \otimes
\overline{h_i h_{i+1}}\otimes...\otimes m =0.
\end{align*}
The second equality follows since we are taking tensor products
over $H_t$. For $i=1, n$ the proof is similar. Hence $\partial_n$
is well-defined, and it is clearly a $H$-module homomorphism. The
proof of ii) is standard and iii) follows by a straightforward
calculation.
\end{proof}

Lemma \ref{tec3} says that the complex
$$\hspace{0.7cm} \dots B_n(H,N)\stackrel{\partial_n}{\longrightarrow
}B_{n-1}(H,N)\dots\longrightarrow
B_2(H,N)\stackrel{\partial_2}{\longrightarrow
}B_{1}(H,N)\stackrel{\partial_1}{\longrightarrow
}B_{0}(H,N)\stackrel{\epsilon}{\longrightarrow } N$$ is acyclic.
Thus, we have a projective resolution of the $H$-module $N$ and we
can compute the Ext groups $\Ext^n_H(N,M)$ for any $M\in \moh$, as
the cohomology groups
$\Ext^n_H(N,M):=\Ker(\partial^n)/\Im(\partial^{n-1})$ of the
complex
$$0\longrightarrow C^0(N,M)\stackrel{\partial^0}{\longrightarrow }C^1(N,M)
\stackrel{\partial^2}{\longrightarrow }\dots\longrightarrow
C^n(N,M)\stackrel{\partial^n}{\longrightarrow
}C^{n+1}(N,M)\longrightarrow\dots$$ where
$C^n(N,M):=\Hom_H(B_n(H,N),M)$.

\medbreak
\subsection{Groupoids }

\

Recall that a (finite) {\it groupoid} is a small category (with
finitely many arrows), such that every morphism has an inverse. We
shall denote a groupoid by $\tgt, \sou:{\mathcal
G}\rightrightarrows {\mathcal P}$, or simply by $\G$, where $\G$
is the set of arrows, $\Pc$ is the set of objects and $\tgt, \sou$
are the target and source maps. The set of arrows between two
objects $P$ and $Q$ is denoted by $\G(P,Q)$ and we shall also
denote $\G(P):=\G(P,P).$ The composition map is denoted by
$m:\G\ftimes{\mathcal G}\rightarrow \G$, and for two composable
arrows $g$ and $h$, that is $\tgt(g)=\sou(h)$, the composition
will be denoted by juxtaposition: $m(g,h)=gh$.

\medbreak

A {\it morphism} between two groupoids is a functor of the
underlying categories. If ${\mathcal G}\rightrightarrows {\mathcal
P}$, and $\Hc\rightrightarrows {\mathcal Q}$ are groupoids, and
$\phi:{\mathcal G}\rightarrow {\mathcal H}$ is a morphism of
groupoids, then for any $q\in {\mathcal P}$,
$\phi(\id_Q)=\id_{\tgt(\phi(\id_Q))}$, and $\phi$ induces a map
$\phi^0:{\mathcal P}\rightarrow {\mathcal Q}$, namely
$\phi^0(p):=\tgt(\phi(\id_Q))$. Thus $\tgt(\phi(g)) =
\phi^0(\tgt(g))$, for any $g\in \G$. If $\phi$ and $\psi$ are two
morphisms of groupoids then $(\phi\psi)^0=\phi^0\psi^0$.

\medbreak We recall a well-known definition.

\begin{defi}Two morphisms of groupoids
$\phi,\psi: {\mathcal G}\rightarrow {\mathcal H}$ are {\it
similar}, denoted $\phi\sim\psi$, if there is a natural
transformation between them; that is, if there exists a map
$\tau:{\mathcal P}\rightarrow {\mathcal H}$ such that
$$\phi(g)\tau(\tgt(g))=\tau(\sou(g))\psi(g), \qquad g\in {\mathcal G}.$$

Observe that "similarity" is an equivalence relation since every
natural transformation between two groupoid morphisms is
necessarily a natural isomorphism.

\medbreak

Two groupoids ${\mathcal G}$, ${\mathcal H}$ are {\it similar},
and we write ${\mathcal G} \sim {\mathcal H}$, if there is an
equivalence of categories between them. In other words, if there
are morphisms $\phi: {\mathcal G}\rightarrow {\mathcal H}$, $\psi:
{\mathcal H}\rightarrow {\mathcal G}$ such that $\phi\circ\psi$
and $\psi\circ\phi$ are similar to the corresponding identities.
\end{defi}

\medbreak

A basic operation between groupoids is the disjoint union. Namely,
if ${\mathcal G}\rightrightarrows {\mathcal P}$, ${\mathcal
H}\rightrightarrows {\mathcal Q}$ are two groupoids, the disjoint
union is the groupoid whose set of arrows is  ${\mathcal G}
\coprod {\mathcal H}$, and whose base is the disjoint union of the
bases: ${\mathcal P} \coprod {\mathcal Q}$. If ${\mathcal G} \sim
{\mathcal G}'$ and ${\mathcal H} \sim {\mathcal H'}$ then
${\mathcal G} \coprod {\mathcal G}' \sim {\mathcal H} \coprod
{\mathcal H}'$.

\medbreak Let us define an equivalence relation on the base $\Pc$
by $P \approx Q$ if $\G(P,Q) \neq \emptyset$. A groupoid $\tgt,
\sou :{\mathcal G}\rightrightarrows {\mathcal P}$ is  {\it
connected } if $P\approx Q$ for all $P, Q\in{\mathcal P}$.

\medbreak

Let $S$ be an equivalence class in $\Pc$ and let $\G_S$ denote the
corresponding connected groupoid with base $S$; that is,
$\G_S(P,Q)=\G(P,Q)$ for any $P,Q\in S$. Then the groupoid
$\mathcal G$ is similar to the disjoint union of the groupoids
$\G_S$: $\G \sim \coprod_{S \in \Pc/\approx} \G_S$.

\medbreak A subgroupoid $\Hc$ of a groupoid $\G$ is \emph{wide} if
$\Hc$ has the same base $\Pc$ as $\G$.

\medbreak

\begin{lema}\label{conexo1} Let $\mathcal G$ be a groupoid.

(i). If $\mathcal G$ is connected, then ${\mathcal G}\sim
{\mathcal G}(P)$ for any $P\in\Pc$.

(ii). If ${\mathcal S}$ is a system of representatives of $\Pc/
\approx$ then
\begin{align}\label{conexo} \G \sim \coprod_{P \in {\mathcal S} }
\G(P).
\end{align}
\end{lema}

\begin{proof} (i). Let us fix $P\in {\mathcal P}$. For any $Q\in
{\mathcal P}$, denote by $\tau_Q$ an element in ${\mathcal
G}(P,Q)$ such that $\tau_{P}=\id_{P}$. So, we have defined a map
$\tau:{\mathcal P}\rightarrow {\mathcal G}$. Define the following
maps $\phi:{\mathcal G}\rightarrow {\mathcal G}(P),$
$\psi:{\mathcal G}(p)\rightarrow {\mathcal G}$ by
$\phi(g)=\tau_{\sou(g)} \; g \; \tau^{-1}_{\tgt(g)}$, $\psi(h)=h$.
In fact these maps are morphisms of groupoids. Since we have
required that $\tau_{P}=\id_{P}$ then $\phi\circ\psi$ is the
identity map. By the definition of $\phi$ we have that
$(\psi\circ\phi)(g)\tau_{\tgt(g)}=\tau_{\sou(g)}g,$ and then
$\psi\circ\phi\sim\id.$ Part (ii) follows from (i).
\end{proof}

\begin{defi} Given a groupoid $\G$ and a map $p:{\mathcal E}\rightarrow\Pc$,
a {\it left action} of $\G$ on $p$ is a map $\fde:\G
{\,}_{\tgt}\times_p \Ec \to \Ec$ such that
\begin{equation*} p(g\fde x)=\sou(g),\qquad
 g \fde(h \fde x)=gh \fde x,\qquad\id_{p(x)} \, \fde  x =  x,
\end{equation*}
for all composable $g,h\in\G$, $x\in\Ec$. We shall say in this
case that $(\Ec, p)$, or $\Ec$, is a \emph{ $\G$-bundle}.

A {\it right action} of $\G$ on ${\mathcal E}$ is a map $\gde :\Ec
{\,}_p\times_{\sou} \G \to \Ec$ such that
\begin{equation*}p(x\gde g)=\tgt(g),\qquad
(x \fde g) \fde h=x \gde gh,\qquad x \gde \id_{p(x)}  =  x,
\end{equation*}
for all composable $g,h\in\G$, $x\in\Ec$.
\end{defi}

\medbreak
\subsection{The groupoid algebra}

\

The {\it groupoid algebra} $\kk\G$ is the $\kk$-algebra with basis
$\{ g: g\in\G\}$,  the product of two elements in the basis being
equal to their composition if they are composable, and $0$
otherwise. The groupoid algebra $\kk\G$ has a weak Hopf algebra
structure via: $\Delta(g)=g\otimes g, \, \epsilon(g)=1,\,
S(g)=g^{-1},$ for all $g\in\G$. The target subalgebra of this weak
Hopf algebra is $\kk \Pc := \oplus_{P\in \Pc} \, \kk \, \id_{P}$.

\medbreak A $\G$-module bundle is a $\G$-bundle $(\Ec,p)$ such
that $\Ec_Q :=p^{-1}(Q)$ is a $\kk$-module for any $Q\in \Pc$ and
the map $g \fde \;:\Ec_{\tgt(g)}\rightarrow \Ec_{\sou(g)}$ is a
linear isomorphism for any $g\in\G$.

\medbreak

There is an equivalence of categories between the category of
${\mathcal G}$-module bundles and $\mog$.

\medbreak The left $\kk\G$-module associated to a $\G$-module
bundle $({\mathcal E},p)$ is given by $M:=\bigoplus_{q\in \Pc}\,
\Ec_Q$, and  the action of $\G$ on $M$ is given by $g.m=g\fde m$
if $m\in \Ec_{\tgt(g)}$ and $g.m=0$ otherwise.  Note that the
fiber $\Ec_Q$ might be zero for some $Q\in \Pc$.

\medbreak Reciprocally, let $M$ be a left $\kk\G$-module and set
$M_P=\id_P M$, $P\in \Pc$; then $M=\bigoplus_{P\in\Pc} M_P$. Let
$$\Ec:=\{(Q,m)\in \Pc\times M \;\text{ such that }\; m\in M_Q\},$$
let $p:\Ec\rightarrow \Pc$ be given by $p(Q,m)=Q$, and let
$\fde:\G {\,}_{\tgt}\times_p \Ec \to \Ec$ be defined by $g\fde
(\tgt(g),m)=(\sou(g),g.m)$. Then $({\mathcal E},p)$ is a
$\G$-module bundle.

\begin{prop}\label{morita} If $\G\rightrightarrows
{\mathcal P}$ and ${\mathcal H}\rightrightarrows {\mathcal Q}$ are
similar groupoids then the categories $\mog$ and $\moha$ are
tensor equivalent. In particular the groupoid algebras $\kk\G$,
$\kk\Hc$ are Morita equivalent.
\end{prop}

\begin{proof}  By hypothesis, there are morphisms of groupoids
$\phi:\G\rightarrow \Hc$ and $\psi:\Hc \rightarrow \G$ satisfying
that $\phi\psi\sim\id_{\Hc}$ and $\psi\phi\sim\id_{\G}$; that is,
there are maps $\theta: \Pc\rightarrow \G$ and
$\eta:\Qc\rightarrow \Hc$ such that
\begin{equation}\label{sim1} \psi\phi(g)\theta(\tgt(g))=\theta(\sou(g))
g,
\end{equation}
\begin{equation}\label{sim2} \phi\psi(h)\eta(\tgt(h))=\eta(\sou(h)) h,
\end{equation}
for all $g\in\G$ and $h\in \Hc$.

\medbreak We define functors $\Psi:\mog\rightarrow \moha$,
$\Phi:\moha\rightarrow \mog$ by
$$\Psi(M)_Q:= M_{\psi^0(Q)},\;\;\;\;\;\text{ and }\;\;\; \Phi(V)_P:=
V_{\phi^0(P),} $$ for all objects $M\in\mog$, $V\in\moha$ and for
all $P\in\Pc$, $Q\in\Qc$. The action of $\Hc$ in $\Psi(M)$ is
defined as follows: if  $x\in\Hc$ and $m\in \Psi(M)_Q$ then
$$x.m=\begin{cases}  \psi(x)m \quad &\text{if } \tgt(x)=Q,\\
  0 \quad &\text{otherwise. }
\end{cases}$$
The action of $\G$ in $\Phi(V)$ is defined as follows: if $g\in\G$
and $v\in\Phi(V)_P$ then
$$g.v=\begin{cases}  \phi(g)v \quad &\text{if } \tgt(v)=P,\\
  0 \quad &\text{otherwise. }
\end{cases}$$
Clearly, these are indeed actions of the corresponding groupoids.
We define natural isomorphisms $\xi:\Id\rightarrow \Psi\Phi$ and
$\zeta:\Psi\Phi\rightarrow\Id $ by
$$\xi_{V_Q}:V_Q\rightarrow V_{(\psi\phi)^0(Q)}, \;\;\;
\xi_{V_Q}(v)=\theta(Q) v \;\;\text{ and}
\;\;\zeta_{M_P}:M_P\rightarrow
M_{(\phi\psi)^0(P)},\;\;\;\zeta_{M_P}(m)=\eta(P)m
$$ for all  $V\in\moha$ and
$M\in\mog$. Equations \eqref{sim1}, \eqref{sim2} imply that $\xi$
and $\zeta$ are morphisms; thus the functors $\Phi$ and $\Psi$
define an equivalence between $\mog$ and $\moha$. A
straightforward verification shows that these functors are in fact
strict tensor functors.
\end{proof}

\medbreak
Combining Lemma \ref{conexo1} and Proposition \ref{morita} we get
\begin{cor} If $\G$ is a connected groupoid then the tensor categories
$\mog$ and $\mogp$ are tensor equivalent for any $P\in \Pc$. In particular the
groupoid algebra $\kk\G$ is Morita equivalent to the group algebra
$\kk \, \G(P)$. \qed
\end{cor}

\medbreak
\subsection{Groupoid cohomology}

\

We briefly recall the well-known groupoid cohomology.

\medbreak Let us fix a groupoid $\;\tgt, \sou:\G \rightrightarrows
\Pc$. Define $\G^{(0)}:=\{\id_Q\}_{Q\in\Pc},\;\; $ $\G^{(1)}=\G$,
and for $n\ge 2$
$$\G^{(n)} =\{(g_1, \dots, g_{n}) \in \G^{n}:
g_1 \vert g_2 \vert \dots \vert g_{n-1}\vert g_{n} \}.$$

Let $(\Ec,p)$ be a $\G$-module bundle, and define
\begin{align*}C^0(\G, \Ec) &=\{f:\Pc\to \Ec:  p(f(Q))=Q \quad \forall\, Q\in
\Pc  \},
\\ C^n (\G, \Ec) &= \{f: \G^{(n)} \to \Ec: f(g_1, \dots, g_{n}) = 0,
\text{ if some }g_i\in \G^{(0)},
\\ & \qquad\qquad\qquad\qquad \text{ and } p(f(g_1, \dots,
g_{n}))=\sou(g_1) \quad \forall\, (g_1, \dots, g_{n}) \in
\G^{(n)}\}.
\end{align*}

The cohomology groups $H^n(\G, \Ec)$ of $\G$ with coefficients in
the $\G$-module bundle $(\Ec,p)$ are the cohomology groups of the
complex

\begin{equation}\label{complex}
0 \longrightarrow C^0 (\G, \Ec) \stackrel{d_0}{\longrightarrow }
C^1 (\G, \Ec) \stackrel{d_1}{\longrightarrow } C^2 (\G,
\Ec)\stackrel{d_2}{\longrightarrow } \dots C^n (\G, \Ec)
\stackrel{d_n}{\longrightarrow } C^{n+1} (\G, \Ec) \longrightarrow
\dots
\end{equation}
where
\begin{equation}\label{cob}
\begin{aligned}d^0 f(g) &= g\fde f(\tgt(g)) - f(\sou(g)), \\
d^n f(g_0, \dots, g_{n}) &=  g_0\fde  f(g_1, \dots, g_{n}) +
\sum_{i=1}^n (-1)^{i} f(g_0, \dots, g_{i-1}g_i, \dots, g_{n})
\\ & \qquad + (-1)^{n+1} f(g_0, \dots, g_{n-1}).
\end{aligned}
\end{equation}

\medbreak Let us denote as usual $Z^n(\G,M):=\Ker(d^n)$,
$B^n(\G,M):=\Im(d^{n-1})$, $n\ge 0$. We next show that this
groupoid cohomology coincides with the cohomology of the weak Hopf
algebra $\kk\G$.

\begin{prop}\label{coho} If $\Ec$ is a $\G$-module bundle and $M$ is the
associated $\kk\G$-module, then the groups $H^n(\G, \Ec)$ and
$H^n(\kk\G, M)$ are naturally isomorphic.

\end{prop}

\begin{proof} Let $C^n(\kk\,\Pc, M):=\Hom_{\kk\G}(B_n(\kk\G,\kk\,\Pc), M)$ as
in Section \ref{weak}. Let us define $F_n:C^n(\kk\,\Pc,
M)\rightarrow C^n(\G,\Ec)$ by
\begin{align*}
F_0(f)(P)  &=  f(\id_P), \qquad \qquad \qquad \qquad \qquad \qquad
\quad P\in \Pc,\\
F_n(f)(g_1,\dots,g_n) &= f(\id_{\sou(g_1)}\otimes
\overline{g_1}\otimes\dots\otimes \overline{g_n} ), \quad
(g_1,\dots, g_n) \in\G^{(n)}.
\end{align*}
Then $F_n$ are isomorphisms whose inverses are the maps
$G_n:C^n(\G,\Ec)\rightarrow C^n(\kk\,\Pc, M)$ given by
\begin{align*}
G_0(f)(g)&=g\fde f(\tgt(g)), \\
G_n(f)(g_o\otimes \overline{g_1}\otimes\dots\otimes
\overline{g_n})&=g_0\fde f(g_1,\dots,g_n).
\end{align*}
It follows from the definition of the maps $d^n$ that $ d^n F_n=
F_n \partial^n$ for any $n \ge 0$. Thus, the maps $F_n$ induce
isomorphisms $H^n(\G, \Ec) \to H^n(\kk\G, M)$.
\end{proof}

As a consequence, we show that groupoid cohomology can be
derived from group cohomology.

\begin{prop} \label{coho2} (i). Let $\Ec$ be a $\G$-module bundle and let
$S$ be a complete set of representatives of equivalence classes in
$\Pc$. Then there are natural isomorphisms-- induced by the
respective inclusions
$$H^n(\G, \Ec)\;\simeq\; \bigoplus_{P\in S }
H^n(\G(P), \Ec_{P}).$$

(ii). Assume that $\G$ is connected and let us fix $O\in \Pc$. Let
$\Ec$ be a $\G$-module bundle and let $\Hc$ be a connected wide
subgroupoid of $\G$. Set $G = \G(O)$, $H = \Hc(O)$. Then the
following diagram commutes:
\begin{align*}
\begin{CD}
H^n(\G, \mathcal{ E}) @>\text{res} >> H^n(\Hc, \Res^{\G}_{\Hc}(\mathcal{ E}))\\
@VVV    @VVV   \\ H^n(G, \mathcal{ E}_O) @>\text{res}>>  H^n(H,
\Res^{G}_{H}( \mathcal{ E}_O)) ,
\end{CD}
\end{align*}
where the vertical arrows are the isomorphisms from part (i).
\end{prop}
\begin{proof} (i). Combine \eqref{conexo}, Proposition
\ref{coho}, and the fact that the $\Ext$ groups are Morita
invariant. Here the well-known natural isomorphisms $\Ext^n_{A
\times B}(M\oplus N, U\oplus V)\simeq \Ext^n_{A}(M, U) \oplus
\Ext^n_{B}(N, V)$, $n\in \Na$, are present, where $A$ and $B$ are
rings, $M$ and $U$ are $A$-modules and $N$ and $V$ are
$B$-modules.

(ii). Straightforward.
\end{proof}

\begin{defi}\label{abarra} Let $A$ be a $\kk$-module. We shall denote by $\underline{A}$
the $\G$-module bundle such that $\underline{A}_P:= A$, $P\in\Pc$,
with trivial action of $\G$. That is, $\underline{A}$ is the
$\G$-module bundle corresponding  to the $\kk\G$-module $\kk\,\Pc
\otimes_\kk A$. By Proposition \ref{coho2} (i), $H^n(\G,
\underline{A})\;\simeq\; \bigoplus_{P\in S} H^n(\G(P), A)$.
\end{defi}

Observe that if $A$ is a $\kk$-module then the set
$Z^2(\G,\underline{A})$ is identified with the set of maps
$\sigma:\G_{\,\tgt}\times_{\sou}\G\to A$ such that
$$ \sigma(g,hf)\sigma(h,f)=\sigma(gh,f)\sigma(g,h), $$
for composable $g,h,f\in\G$.

\bigbreak
\section{The Kac exact sequence for matched pairs of groupoids }

\subsection{Matched pairs of groupoids}

\

We briefly recall the definition of matched pair of groupoids, and
the equivalent formulations in terms of exact factorizations or
vacant double groupoids, see \cite{Ma} or \cite{AN} for details.

\medbreak A {\it matched pair of groupoids} is a collection
$(\Hc,\Vc,\fde,\gde)$, where  $b, t:\Vc\rightrightarrows \Pc$ and
$r,l:{\mathcal H}\rightrightarrows \Pc$ are two groupoids over the
same base $\Pc$, $\fde:{\mathcal H}\hvtimes \Vc\rightarrow
{\mathcal V}$ is a left action of $\Hc$ on $(\Vc, t)$,
$\gde:{\mathcal H}\hvtimes \Vc\rightarrow \Hc$ is a right action
of $\Vc$ on $(\Hc, r)$ such that

\begin{equation} \label{m7}  b(x\fde g)= l(x\gde g),\quad
x\fde gh =(x\fde g)((x\gde g)\fde h),\quad
 xy\gde g = (x\gde (y\fde g))(y\gde g),
\end{equation}
for composable elements $x,y\in\Hc$ and $g,h\in\Vc$. Here and
below we use the `horizontal and vertical notation": the source
and target of $\Hc$, resp. $\Vc$, are denoted $l$ and $r$ (left
and right), resp. $t$ and $b$ (top and bottom).

\medbreak Let $(\Hc,{\mathcal V}, \gde, \fde)$ be a matched pair
of groupoids. There is an associated {\it diagonal groupoid} $\Do$
with set of arrows $\Vc\dtimes \Hc$, base $\Pc$, source, target,
composition and identity given by
\begin{equation*}
\sou(g,x)= t(g),\quad  \tgt(g,x)= r(x), \quad (g,x)(h,y)=(g(x\fde
h),(x\gde h)y), \quad \id_P =(\id_P,\id_P),
\end{equation*}
$g,h\in \Vc$, $x,y\in \Hc$, $P\in \Pc$. Then we have an exact
factorization of groupoids $\Do = \Vc\Hc$. Conversely, if $\D =
\Vc\Hc$ is an exact factorization of groupoids then there are
actions $\gde$, $\fde$ such that $(\Hc,{\mathcal V}, \gde, \fde)$
form a matched pair of groupoids, and $\D \simeq \Do$.

\medbreak There is also a vacant double groupoid associated to the
matched pair of groupoids $(\Hc, \Vc)$. In simple terms, this is a
collection of groupoids $\begin{matrix} {\mathcal B}
&\rightrightarrows &{\mathcal H}
\\\downdownarrows &&\downdownarrows \\ {\mathcal V} &\rightrightarrows &\Pc
\end{matrix}$
with the following meaning. A pair $(x,g)$ in ${\mathcal B}:=\Hc
\hvtimes\Vc$ is depicted as a box $\begin{matrix}   x \\ h\,\,
\boxe \,\,g   \\ y
\end{matrix}$
where $h=x\fde g$, $y=x\gde g$.  The horizontal groupoid
$r,l:{\mathcal B}\rightrightarrows \Vc$ has source, target,
composition and identity given by
$$r(x,g)=g,\quad l(x,g)=x\fde g,\quad (x,g)(y,h)=(xy,h), \quad
\id(g)=(\id_{t(g)},g) $$ $x, y\in\Hc$, $g, h\in\Vc$. The vertical
groupoid $t,b:{\mathcal B} \rightrightarrows \Hc$ has source,
target, composition and identity given by
$$b(x,g)=x\gde g,\quad t(x,g)=x,\quad (x,g)(y,h)=(x,gh),\quad \id(x)
=(x,\id_{r(x)}),$$ $x, y\in\Hc$, $g, h\in\Vc$. If $A, B$ are two
boxes we denote $A\mid B$ if they are horizontally composable, and
$ \displaystyle \frac{A}{B}$ if they
are vertically composable, that is if $A=\begin{matrix}   x \\
h\,\,  \boxe \,\,g   \\ y\end{matrix}$ and $B=\begin{matrix}   z
\\ f\,\, \boxe \,\,k   \\ w
\end{matrix}$, then $A\mid B$ if and only if $g=f$, and
$ \displaystyle \frac{A}{B}$ if and only if $y=z$.

\medbreak
\subsection{The Kac Exact Sequence}

\

In this Subsection we shall review and complete details of the
proof of the Kac exact sequence for vacant double groupoids
introduced in \cite{AN}. Let $(\Hc,{\mathcal V}, \gde, \fde)$ be a
matched pair of groupoids. We begin by  a non standard resolution
of the diagonal groupoid, adapting ideas from \cite{M1} to the
groupoid case. If $r,s \in \Na$, we denote by $\B^{[r,s]}$ the set
of matrices
$$
\left(\begin{tabular}{p{0,8cm} p{0,8cm} p{0,8cm} p{0,8cm}}
$A_{11}$ & $ A_{12}$ & \dots & $A_{1s}$ \\
$A_{21}$ & $ A_{22}$ & \dots & $A_{2s}$ \\
\dots & \dots & \dots & \dots \\
$A_{r1}$ & $ A_{r2}$ & \dots & $A_{rs}$  \end{tabular}\right) \in
\B^{r\times s}
$$
such that
\begin{itemize}
\item For all $i,j$, $A_{ij}\mid A_{i,j+1}$, $ \displaystyle
\frac{A_{ij}}{A_{i+1,j}}$. This condition is summarized in the
notation:

\centerline{
\begin{tabular}{p{0,8cm}|p{0,8cm}|p{0,8cm}|p{0,8cm}}
$A_{11}$ & $ A_{12}$ & \dots & $A_{1s}$ \\ \hline $A_{21}$ & $
A_{22}$ & \dots & $A_{2s}$ \\ \hline \dots & \dots & \dots & \dots
\\ \hline $A_{r1}$ & $ A_{r2}$ & \dots & $A_{rs}$.  \end{tabular}}

\smallbreak

\item  If $j< s$ then $A_{ij}$ is not a horizontal identity.
\item If $i > 1$ then $A_{ij}$ is not a vertical identity.
\end{itemize}

\medbreak Observe that if $A=(A_{ij})$ is an element in
$\B^{[r,s]}$ then $A$ is determined by $s$ composable elements
$x_1,\dots, x_s$ in $\Hc$--those in the top of the array-- and $r$
composable elements $g_1,\dots, g_r$ in $\Vc$--those in the right
side of the array, with $r(x_s) = t(g_1)$. We shall denote it by
$$A=: \begin{matrix} x_1,\dots,x_s &
\\ \begin{tabular}{|p{1,4cm}|}
\hline \\ \hline \end{tabular} & g_1,\dots,g_r \\
\\\end{matrix}.$$ Let $s^{r,s}_V:\kk\B^{[r,s]}\rightarrow
\kk\B^{[r+1,s]}$ (vertical homotopy maps),
$\partial^{r,s}_V:\kk\B^{[r,s]}\rightarrow \kk\B^{[r-1,s]}$ and
$\partial^{r,s}_H:\kk\B^{[r,s]}\rightarrow \kk\B^{[r,s-1]}$
(vertical and horizontal coboundary maps) be defined by

\begin{align*}
s^{r,s}_V\left(\begin{matrix} x_1,\dots,x_s &
\\ \begin{tabular}{|p{1,4cm}|}
\hline \\ \hline \end{tabular} & g_1,\dots,g_r \\ \\
\end{matrix}\right) &:= \begin{matrix} x_1,\dots,x_s
& \\
\begin{tabular}{|p{1,4cm}|}
\hline \\ \hline \end{tabular} & \id_{t(g_1)},g_1,\dots,g_r \\ \\
\end{matrix};
\end{align*}
\begin{align*}
\partial_H^{r,s}\left(\begin{matrix} x_1,\dots,x_s &
\\ \begin{tabular}{|p{1,4cm}|}
\hline \\ \hline \end{tabular} & g_1,\dots,g_r \\ \\
\end{matrix}\right)&=
  \sum_{j=1}^{s-1}(-1)^{s-j-1} \, \begin{matrix} x_1,..,x_jx_{j+1},..,x_s
&\\ \begin{tabular}{|p{2,7cm}|}
\hline \\ \hline \end{tabular} & g_1,\dots,g_r \\ \\
\end{matrix}\\
& \quad + (-1)^{s-1} \, \begin{matrix} x_2,\dots,x_s
& \\
\begin{tabular}{|p{1,4cm}|}
\hline \\ \hline \end{tabular} & g_1,\dots,g_r\\  \\
\end{matrix};
\end{align*}
\begin{align*}
\partial_V^{r,s}\left(\begin{matrix}x_1,\dots,x_s &
\\ \begin{tabular}{|p{1,4cm}|}
\hline \\ \hline \end{tabular} & g_1,\dots,g_r \\ \\
\end{matrix}\right)&=\sum_{i=1}^{r-1}(-1)^{i-1} \,
\begin{matrix} x_1,\dots,x_s
& \\
\begin{tabular}{|p{1,4cm}|} \hline \\ \hline \end{tabular}
& g_1,\dots,g_ig_{i+1},\dots,g_r \\ \\
\end{matrix}\\
&+ (-1)^{r-1} \, \begin{matrix}x_1,\dots,x_s
&\\
\begin{tabular}{|p{1,4cm}|} \hline \\ \hline \end{tabular}
& g_1,\dots,g_{r-1}\\ \\\end{matrix}.
\end{align*}

A straightforward computation shows that the following diagram
commutes:
$$
\begin{CD}
\kk\B^{[r, s]} @>{\partial^{r,s}_H}>> \kk\B^{[r, s-1]}
\\
@V{\partial^{r,s}_V}VV @VV{\partial^{r,s-1}_V}V \\
\kk\B^{[r-1, s]} @>{\partial^{r-1,s}_H}>> \kk\B^{[r-1, s-1]}.
\end{CD}
$$
Thus, we have constructed a double chain complex $\Bg^{\bullet,
\bullet}$:
$$ \Bg^{\bullet, \bullet} = \qquad
\begin{CD}
\vdots \\
@VVV \\
\kk\B^{[3, 1]} @<{\partial_H}<< \hspace{18pt}\raisebox{1.0ex}{\vdots}\cdots \\
@VV{\partial_V}V @VV{-\partial_V}V \\
\kk\B^{[2, 1]} @<{\partial_H}<< \kk\B^{[2, 2]} @<{\partial_H}<<
\hspace{18pt}\raisebox{1.0ex}{\vdots}\cdots \\
@VV{\partial_V}V @VV{-\partial_V}V @VVV \vspace{9pt}\\
\kk\B^{[1, 1]} @<{\partial_H}<< \kk\B^{[1, 2]} @<{\partial_H}<<
\kk\B^{[1, 3]} @<<< \cdots\ \
\end{CD}
$$

\medbreak Note that $\B^{[r,s]}$ is a subset of the set
$\B^{(r,s)}$ defined in \cite{AN}; and the double complex
presented here is different from the double complex in \cite{AN}.
However these will be reconciliated in Remark \ref{comparacion}
below.

Let us now define an action of the diagonal groupoid $\Do$ on
$\kk\B^{[r,s]}$, $r, s
> 0$. The $\Pc$-gradation on $\kk\B^{[r,s]}$ is given by:
$$p \left(\begin{matrix} x_1,\dots,x_s &
\\ \begin{tabular}{|p{1,4cm}|}
\hline \\ \hline \end{tabular} & g_1,\dots,g_r \\
\\\end{matrix}\right):= t(g_1)= r(x_s).$$

If $h\in\Vc$, $y\in\Hc$ are such that $r(y)= r(x_s)= t(g_1) =
b(h)$, then we set

\begin{align}\label{accion1}
y.\left(\begin{matrix} x_1,\dots,x_s &
\\ \begin{tabular}{|p{1,4cm}|}
\hline \\ \hline \end{tabular} & g_1,\dots,g_r \\ \\\end{matrix}
\right)&:=\begin{matrix} x_1,\dots,x_sy^{-1} &
\\
\begin{tabular}{|p{1,9cm}|}
\hline \\ \hline \end{tabular} & y\fde g_1,
(y\gde g_1)\fde g_2,\dots,(y\gde g_1\dots g_{r-1})\fde g_r \\ \\
\end{matrix},
\\
\label{accion2} h . \left(\begin{matrix}\begin{matrix}
x_1,\dots,x_s &
\\ \begin{tabular}{|p{1,4cm}|}
\hline \\ \hline \end{tabular} & g_1,\dots,g_r \\ \\\end{matrix}
\end{matrix}\right) &:=
\begin{matrix} x_1\gde (x_2\dots x_s\fde h^{-1}),\dots,x_s\gde h^{-1}
&\\
\begin{tabular}{|p{4,9cm}|}
\hline \\ \hline \end{tabular} & hg_1,g_2,\dots,g_r \\ \\
\end{matrix}.
\end{align}

\begin{lema}\label{tec2} \begin{itemize}
\item[(i)] The rules \eqref{accion1} and \eqref{accion2} induce a
structure of $\kk\Do$-module on $\kk\B^{[r,s]}$.

\medbreak \item[(ii)] There are $\kk\Do$-isomorphisms:

\begin{itemize}
\item[$\bullet $] $\kk\B^{[r,s]}\simeq
\kk\Do\otimes_{\kk\Pc}\kk\B^{[r-1, s-1]} \text{ for any }\; r,s >
1,$ \item[$\bullet $] $\kk\B^{[r, 1]}\simeq \kk\Do\otimes_{\kk\Pc}
\kk\Vc^{(r-1)}$ for any $r> 0$, \item[$\bullet $] $ \kk\B^{[1,
s]}\simeq \kk\Do\otimes_{\kk\Pc} \kk\Hc^{(s-1)}$ for any $s> 0$,
\end{itemize}
\end{itemize}

\smallbreak

\noindent
where the action of $\kk\Do$ on $\kk\Do\otimes_{\kk\Pc} L$ is on
the first tensorand, for any  $L \in \,_{\kk\Pc}{\mathcal M}$.

\medbreak \begin{itemize}
\item[(iii)] $\kk\B^{[r,s]}$ is a projective
$\kk\Do$-module for any $r,s > 0$.

\medbreak \item[(iv)] The coboundary maps $\partial^{r, s}_V,
\partial^{r, s}_H$ are morphisms of $\kk\Do$-modules.

\medbreak \item[(v)] $\partial^{r+1,s}_V s^{r,s}_V +
s^{r-1,s}_V\partial^{r,s}_V= \id_{\kk\B^{[r,s]}}.$
\end{itemize}
\end{lema}

\begin{proof} (i) Let $A = \begin{matrix} x_1,\dots,x_s &
\\ \begin{tabular}{|p{1,4cm}|}
\hline \\ \hline \end{tabular} & g_1,\dots,g_r \\ \\\end{matrix}
\in \B^{[r,s]}$; let $h\in \Vc, y\in \Hc $, such that $r(y)=t(h)$
and $b(h)=p(A)$. We claim that  $y.(h.A)= (y\fde h). \left((y\gde
h). A\right)$. We have
$$y.(h.A)=
\begin{matrix} x_1\gde(x_2\dots x_s\fde h^{-1}),\dots,(x_s\gde
h^{-1})y^{-1}
&\\
\begin{tabular}{|p{5,5cm}|}
\hline \\ \hline \end{tabular} & y\fde hg_1,
(y\gde hg_1)\fde g_2,\dots,(y\gde hg_1\dots g_{r-1})\fde g_r \\ \\
\end{matrix};$$
and $(y\fde h). \left((y\gde h). A\right) =$
$$\begin{matrix} x_1\gde((x_2\dots
x_s(y\gde h)^{-1})\fde (y\fde h)^{-1}),\dots,(x_s(y\gde
h)^{-1})\gde (y\fde
h)^{-1} &\\
\begin{tabular}{|p{8cm}|}
\hline \\ \\\hline \end{tabular} & \hspace{-27pt} (y\fde h)((y\gde
h)\fde
g_1),\dots,((y\gde h g_1\dots g_{r-1})\fde g_r. \\
\end{matrix}
$$
Then $y.(h.A)=(y\fde h). \left((y\gde h). A\right)$ by \eqref{m7}
and the identities $$(y\fde h)^{-1}=(y\gde h)\fde h^{-1}, \quad
(y\gde h)^{-1}=y^{-1}\gde (y\fde h).$$

(ii) Assume that $r,s > 1$. Define the maps $\phi:\kk\B^{[r,s]}\to
\kk\Do\otimes_{\kk\Pc}\kk\B^{[r-1,s-1]}$, $\psi:
\kk\Do\otimes_{\kk\Pc}\kk\B^{[r-1,s-1]}\to \kk\B^{[r,s]}$ by the
formulas

\begin{multline*}\phi\left(\begin{matrix} x_1,\dots,x_s
\\ \begin{tabular}{|p{1,4cm}|}
\hline \\ \hline \end{tabular} & g_1,\dots,g_r \\
\\\end{matrix}\right):=(x^{-1}_s,x_s\fde g_1)\,
\otimes\\
\begin{matrix} x_1\gde(x_2...x_s\fde g_1),\dots,x_{s-1}\gde
(x_s\fde g_1)
\\ \begin{tabular}{|p{5,5cm}|}
\hline \\ \hline \end{tabular} & (x_s\gde g_1)\fde g_2,\dots,
(x_s\gde g_1...g_{r-1})\fde g_r \\
\end{matrix}, \end{multline*}

\begin{align*}\psi((y,h)\otimes \begin{matrix} x_1,\dots,x_{s-1} &
\\ \begin{tabular}{|p{1,8cm}|}
\hline \\ \hline \end{tabular} & g_1,\dots,g_{r-1} \\
\\\end{matrix}):= (y,h).\left(\begin{matrix} x_1,\dots,x_{s-1}, \id_{r(x_{s-1})} &
\\ \begin{tabular}{|p{2,9cm}|}
\hline \\ \hline \end{tabular} & \id_{t(g_1)},g_1,\dots,g_{r-1} \\
\\\end{matrix}\right).
\end{align*}
These maps are morphisms of $\kk\Do$-modules and one is the
inverse of each other. The proof of the cases $r=1$ or $s=1$
follows similarly. Part (iii) follows from (ii) and Lemma
\ref{tec}. The proof of (iv) and (v) is straightforward.

\end{proof}

\medbreak Let ${\mathcal A}^{\bullet, \bullet}$ be the double
chain complex obtained from $\Bg^{\bullet, \bullet}$ by removing
the edges; that is ${\mathcal A}^{r,s}:=\Bg^{r+1,s+1}$. Let $M$ be
a left $\kk\Do$-module and thus a left $\kk\Hc$-module and also a
left $\kk\Vc$-module. Define the double cochain complexes
$\Bg^{\bullet, \bullet}(M), {\mathcal E}^{\bullet, \bullet}(M),
{\mathcal A}^{\bullet, \bullet}(M)$ by

$$\Bg^{r,s}(M):= \Hom_{\kk\Do}(\kk\B^{[r,s]}, M),\,\,\,\,\,\,\,
{\mathcal A}^{r,s}(M):=\Hom_{\kk\Do}(\kk\B^{[r+1,s+1]}, M),$$ and
${\mathcal E}^{r,s}(M)$ consists only of the edges of $\Bg(M)$.

\begin{rmk}\label{identificaciones} Let $M$ be a $\kk\Do$-module. Lemma \ref{tec2}
(ii) implies that there are  natural $\kk$-linear isomorphisms:
$\Bg^{r,s}(M)\simeq \Hom_{\kk}(\kk\B^{[r-1,s-1]},M) $ for any $r,s
> 1$, and there are natural bijections: $\Bg^{r,1}(M)\simeq C^{(r-1)}(\Vc,\Ec),\,$
$\Bg^{1,s}(M)\simeq C^{(s-1)}(\Hc,\Ec)$ for any $r,s> 0$, where
$\Ec$ is the module bundle corresponding to $M$.
\end{rmk}

\begin{rmk}\label{normalizacion} Let $\B^{(r,s)}$
be as in \cite{AN}. Suppose that $r,s
> 1$. We extend any $\mu \in \Bg^{r,s}(M) \simeq
\Hom_{\kk}(\kk\B^{[r-1,s-1]},M)$ to $\widetilde \mu \in
\Hom_{\kk}(\kk\B^{(r-1,s-1)},M)$ by 0 on $\B^{(r-1,s-1)} -
\B^{[r-1,s-1]}$. In other words, the elements of $\Bg^{r,s}(M)$
are normalized by definition.
\end{rmk}

\medbreak

Now we can formulate the Kac exact sequence for groupoids.

\begin{teo}\cite[Prop 3.14]{AN} Let $M$ be a $\kk\Do$-module.
Then, there is an exact sequence
\begin{equation}
\begin{aligned}
0 &\longrightarrow H^1(\D, M)\stackrel{\text{ res }}
{\longrightarrow} H^1(\Hc, M) \oplus H^1(\Vc, M)
\longrightarrow H^0(\Tot {\mathcal A}^{\bullet, \bullet}(M)) \\
&\longrightarrow H^2(\D, M)\stackrel{\text{ res }}
{\longrightarrow} H^2(\Hc, M) \oplus H^2(\Vc, M) \longrightarrow
H^1(\Tot {\mathcal A}^{\bullet, \bullet}(M)) \\
&{\longrightarrow} H^3(\D, M)\stackrel{\text{ res
}}{\longrightarrow} H^3(\Hc,M) \oplus H^3(\Vc, M)\longrightarrow
\dots
\end{aligned}
\end{equation}
The maps denoted by $res$ in the above exact sequence come from
the usual restriction maps.
\end{teo}

\begin{proof} The short exact sequence of double complexes
$0\to {\mathcal A}^{\bullet, \bullet}(M)\to \Bg^{\bullet,
\bullet}(M)\to \mathcal{E}^{\bullet, \bullet}(M)\to 0$ induces a
long exact sequence in cohomology. By remark
\ref{identificaciones} it is easy to see that
$$H^n(\Tot\mathcal{E}^{\bullet, \bullet}(M))\simeq
H^n(\Hc,M)\bigoplus H^n(\Vc,M)$$ for any $n\in\Na_0$. By  Lemma
\ref{tec2} (v) each column  of $\Bg^{\bullet,\bullet}$ is acyclic.
Hence the associated total complex
\begin{equation} \longrightarrow\Tot(\Bg)_n
\stackrel{\partial_n}{\longrightarrow
}\Tot(\Bg)_{n-1}\dots\longrightarrow
\Tot(\Bg)_3\stackrel{\partial_1}{\longrightarrow}
\Tot(\Bg)_2\stackrel{\epsilon}{\longrightarrow } \kk\,\Pc,
\end{equation}
where $\epsilon: \Tot(\Bg)_1\rightarrow \kk\,\Pc$ is given by the
degree: $\epsilon\Bigl(\begin{matrix}   x \\  \,\, \boxe \,\, g \\
\quad
\end{matrix}\Bigr)= l(x)$, is a projective resolution of the trivial
$\kk\Do$-module. See, for instance, \cite[Ex. 1.2.5]{W}. Hence
$H^n(\Tot\Bg(M))\simeq H^n(\Do,M)$ for any $n\in\Na_0$, and this
finishes the proof.
\end{proof}

\begin{rmk}\label{comparacion} Let now $\kk=\Z$ and $A = \ku^{\times}$ and recall the
meaning of $\underline{\ku^{\times}}$ in Definition \ref{abarra}.
Let us denote
$$\Aut(\ku\, \Tc):=H^0(\Tot {\mathcal
A}^{\bullet, \bullet}(\underline{\ku^{\times}})), \qquad
\Opext(\Vc,\Hc):=H^1(\Tot {\mathcal A}^{\bullet,
\bullet}(\underline{\ku^{\times}})).$$

The group $Z^2(\Tot {\mathcal A}^{\bullet,
\bullet}(\underline{\ku^{\times}}))$ can be identified with set of
pairs $(\sigma,\tau)$ such that $\sigma$ is a normalized 2-cocycle
with values in $\underline{\ku^{\times}}$ for the vertical
groupoid $\mathcal{B}\rightrightarrows \Hc$, $\tau$ is a
normalized 2-cocycle with values in $\underline{\ku^{\times}}$ for
the horizontal groupoid $\mathcal{B}\rightrightarrows \Vc$ and
\begin{equation}\label{cociclo-sigma-tau}
\quad \sigma(AB, CD) \tau \left(\begin{matrix} A \\ C
\end{matrix},
\begin{matrix}B \\ D \end{matrix}\right)
= \tau(A, B) \tau(C, D) \sigma(A, C) \sigma(B, D),  \end{equation}
for  any $A, B, C, D$ such that $ \quad
\begin{tabular}{p{0,4cm}|p{0,4cm}} A & B
\\ \hline C & D \end{tabular}.$ Hence, $\Opext(\Vc,\Hc)$ coincides
with the group considered in \cite{AN}. We have the familiar
expression
\begin{equation}\label{kac}
\begin{aligned}
0 &\longrightarrow  H^1(\D, \underline{\ku^{\times}})
\stackrel{\text{ res}}{\longrightarrow} H^1(\Hc,
\underline{\ku^{\times}}) \oplus H^1(\Vc,
\underline{\ku^{\times}})
\longrightarrow \Aut(\ku\,\Tc) \\
&\longrightarrow  H^2(\D, \underline{\ku^{\times}})
\stackrel{\text{ res}}{\longrightarrow} H^2(\Hc,
\underline{\ku^{\times}}) \oplus
H^2(\Vc, \underline{\ku^{\times}}) \longrightarrow  \Opext(\Vc,\Hc) \\
&\longrightarrow H^3(\D, \underline{\kk^{\times}})\stackrel{\text{
res }}{\longrightarrow} H^3(\Hc,\underline{\ku^{\times}}) \oplus
H^3(\Vc, \underline{\ku^{\times}})\longrightarrow \dots
\end{aligned}
\end{equation}

Notice that, if $(\sigma, \tau)\in Z^2(\Tot {\mathcal A}^{\bullet,
\bullet} (\underline{\ku^{\times}}))$, then it follows from
equation \eqref{cociclo-sigma-tau} that
\begin{equation}\label{lemaprevio}
\sigma((\id_{\sou(g)},g),(\id_{\sou(h)},h))=1, \quad
\tau((x,\id_{\tgt(x)}),(y,\id_{\tgt(y)}))=1, \qquad g,h\in\Vc,
\quad x,y\in\Hc.
\end{equation}
\end{rmk}

\bigbreak
\section{Matched Pairs of groupoids with Connected Vertical groupoid}

Let $\D\rightrightarrows \Pc$ be a connected groupoid. We fix
$O\in\Pc$ and $\tau_P\in \D(O, P)$, $P\in\Pc$, $\tau_O = \id_O$.
We denote by $D$ the group $\D(O)$; thus $\D\simeq D\times \Pc^2$,
with isomorphism given by
$$
\D(P, Q) \ni x \mapsto (\tau_P x\tau_Q^{-1}, (P, Q)).
$$
In other words, we pull back $x$ to an arrow from $O$ to $O$ via the $\tau$'s:
\[\xymatrix@C+5pt@R-5pt{
{O}\ar@/_/[d]_{}\ar@/^/[r]^{\tau_P} & {P}\ar@/^/[d]^{x}\\
{O}\ar@/^/[r]_{\tau_Q} & {Q} }\] Different choices of families
$\tau_P\in \D(O, P)$, $P\in\Pc$, just amount to different
isomorphisms of groupoids $\D\simeq D\times \Pc^2$.

\medbreak
\subsection{Structure of exact factorizations}

\

Let us fix an equivalence relation $\approx_H$ on $\Pc$ and a
section $\sigma:\Pc/\approx_H\ \to \Pc$ of the canonical
projection. Then there is a bijection between

\begin{itemize}
\item[(a)]
the set of wide subgroupoids $\Hc$ of $\D$ with equivalence
relation $\approx_H$, and

\item[(b)] the set of collections $( (H_{\clase})_{\clase\in  \Pc/\approx_H},
(\overline{\lambda_P})_{P\in \Pc})$, where
 $H_{\clase}$ is a subgroup of $D$, $\clase\in  \Pc/\approx_H$, and
$\overline{\lambda_P}\in H_{\clase}\backslash D$, $P\in \clase$,
with $\lambda_{\sigma(\clase)}\in H_{\clase}$ for any $\clase\in
\Pc/\approx_H$.
\end{itemize}
Namely, from (a) to (b), if $\clase\in  \Pc/\approx_H$ and  $P\in
\clase$, then we choose $g_P\in \Hc(\sigma(\clase), P)$ and set
$$H_{\clase} = \tau_{\sigma(\clase)}\Hc(\sigma(\clase))\tau_{\sigma(\clase)}^{-1},
\qquad \lambda_P = \tau_{\sigma(\clase)}g_P\tau_{P}^{-1}.
$$
The choice of $g_P$ does not affect the class of $\lambda_P$ in
$H_{\clase}\backslash D$. By definition,
$\lambda_{\sigma(\clase)}\in H_{\clase}$.

\medbreak Conversely, if a collection $( (H_{\clase})_{\clase\in
\Pc/\approx_H}, (\overline{\lambda_P})_{P\in \Pc})$ satisfies the
above conditions then the wide subgroupoid $\Hc$ that this
collection determines is given by
$$
\Hc(P,Q)=
\begin{cases} \tau_P^{-1}\lambda^{-1}_P H_{\clase} \lambda_Q \tau_Q, \quad\,
\text{ if } \quad
P\approx_H Q, P, Q\in \clase;\\
\emptyset \qquad\qquad\qquad\qquad \text{ if } \quad
P\not\approx_H Q.
\end{cases}
$$
In other words, \emph{cf.} the equality $\lambda_P\tau_{P} =
\tau_{\sigma(\clase)}g_P$, arrows $x$ in $\Hc$ from $P$ to $Q$
correspond to arrows $\widetilde x \in H_{\clase}$ as in this
commutative diagram:
\[\xymatrix@R-5pt@C+15pt{
{O}\ar@/_/[d]_{\widetilde x}\ar@/^/[r]^{\tau_{\sigma(\clase)}} & {\sigma(\clase)}\ar@/_/[d]_{}\ar@/^/[r]^{g_P} & {P}\ar@/^/[d]^{x}\\
{O}\ar@/^/[r]_{\tau_{\sigma(\clase)}} &
{\sigma(\clase)}\ar@/^/[r]_{g_Q} & {Q} }\]

\medbreak By abuse of notation, we shall say that $\Hc$ is
associated to $((H_{\clase})_{\clase\in  \Pc/\approx_H},
(\overline{\lambda_P})_{p\in \Pc})$.

\medbreak We reformulate the description of exact factorizations
of a connected groupoids given in \cite[Th. 2.15]{AN} in terms of
the preceding discussion. Let us fix equivalence relations
$\approx_V$ and $\approx_H$ on $\Pc$ and sections
$\rho:\Pc/\approx_V \to \Pc$, $\sigma:\Pc/\approx_H \to \Pc$ of
the canonical projections.

\begin{teo}\label{estructdouble}
Let $\Hc$ and $\Vc$ be wide subgroupoids of $\D$ associated
respectively  to collections $((H_{\clase})_{\clase\in
\Pc/\approx_H}, (\lambda_P)_{P\in \Pc})$,
$((V_{\clate})_{\clate\in \Pc/\approx_V}, (\mu_P)_{P\in\Pc})$ as
explained above. Then the  following are equivalent:

\begin{itemize}
\item[(i)] $\D = \Vc \Hc$ is an exact factorization.

\medbreak \item[(ii)] $(\Hc, \Vc)$ is a matched pair of groupoids
and  $\D \simeq \Vc \bowtie\Hc$.

\medbreak \item[(iii)] The  following conditions hold:

\begin{align}\label{tdoubleclasses}
D &= \coprod_{R\in\, \clase \cap \clate} V_{\clate} \, \mu_R
\lambda_R^{-1} \, H_{\clase}, \qquad \text{ for all } \clase\in
\Pc/\approx_H, \clate\in  \Pc/\approx_V;
\\ \label{tintersecciones}
\mu^{-1}_P V_{\clate} \mu_P \cap \lambda^{-1}_P
H_{\clase}\lambda_P &= \{ 1\}, \qquad\qquad\qquad \text{ for all }
\clase\in \Pc/\approx_H,\, \clate\in  \Pc/\approx_V, P\in \clase
\cap \clate.
\end{align}
\qed
\end{itemize}
\end{teo}

\medbreak In the rest of this Section we shall study matched pairs
of groupoids whose vertical groupoid is connected. These are
exactly those such that the category of finite dimensional
representations of the weak Hopf algebra $\ku^{\tau}_{\sigma} \Tc$
is fusion, that is the unit object is simple  \cite[Prop.
3.10]{AN}. See \cite{ENO} and references therein for the notion of
fusion category.

\medbreak We fix a vertical subgroupoid $\Vc$ of $\D$ determined
by a subgroup $V$ of $D$ and $\overline{\mu_P} \in V\backslash D$,
$P\in\Pc$. Without loss of generality, we can assume that
$\mu_P=1$ for all $P\in\Pc$; just change the family $(\tau_P)$ by
$(\mu_P \tau_P)$.

\medbreak We shall say that $\Hc$ is an \emph{exact factor} of the
subgroupoid $\Vc$ of $\D$ if $\D = \Vc \Hc$ is an exact
factorization.

\begin{cor}\label{estructdoubleconn}
Let $\Hc$ be a wide subgroupoid of $\D$ associated to
$((H_{\clase})_{\clase\in  \Pc/\approx_H}, (\lambda_P)_{P\in
\Pc})$. Then the  following are equivalent:

\begin{itemize}
\item[(i)] $\Hc$ is an exact factor of $\Vc$.

\medbreak \item[(ii)] The  following conditions hold:

\begin{align}\label{doubleclasses}
D &= \coprod_{R\in\, \clase} V \,  \lambda_R^{-1} \, H_{\clase},
\qquad \text{ for all } \clase\in   \Pc/\approx_H;
\\ \label{intersecciones}
V  \cap g H_{\clase} g^{-1} &= \{ 1\}, \qquad\qquad\qquad \text{
for all } g\in  D.
\end{align}
\end{itemize}
\end{cor}

\pf Condition \eqref{tintersecciones} reads now $V  \cap
\lambda^{-1}_P H_{\clase}\lambda_P = \{ 1\}$, for all $\clase\in
\Pc/\approx_H$, $P\in s$. In presence of \eqref{doubleclasses},
this is equivalent to \eqref{intersecciones}. \epf

\begin{sm}\label{genconst}
To construct an explicit example of an exact factorization $\D =
\Vc\Hc$ with $\Vc$ connected, we need:

\begin{itemize}
\item A finite group $D$, a subgroup $V$ of $D$
and a finite non-empty set $\Pc$.

\medbreak We fix $O\in \Pc$ and define $\D$ and $\Vc$ as explained
above.

\medbreak \item An equivalence relation $\approx_H$ in $\Pc$.

\medbreak \item A family  $(H_{\clase})_{\clase\in \Pc/\approx_H}$
of subgroups of $D$ such that

\medbreak
\begin{itemize}
\item[(a)] $V$ intersects trivially all conjugates of $H_{\clase}$
for all $\clase\in  \Pc/\approx_H$.

\medbreak
\item[(b)] There are bijections $\varphi_{\clase}: V\backslash D/H_{\clase} \simeq
\clase$ for all $\clase\in  \Pc/\approx_H$.

\medbreak We denote $\sigma(\clase) =
\varphi_{\clase}(VH_{\clase}) \in \clase$.
\end{itemize}

\medbreak \item A section $\zeta:\clase \to D$ of the canonical
projection $D \to V\backslash D/H_{\clase}$ composed with
$\varphi_{\clase}$, such that $\zeta_{\sigma(\clase)} \in
H_{\clase}$.
\end{itemize}

\medbreak We set $\lambda_P = \zeta^{-1}_P$ for any $P\in \Pc$;
and clearly $\lambda_{\sigma(\clase)} =
\zeta^{-1}_{\sigma(\clase)} \in H_{\clase}$. Then $\Hc$ is the
wide subgroupoid associated to $((H_{\clase})_{\clase\in
\Pc/\approx_H}, (\lambda_P)_{P\in \Pc})$; it is an exact factor of
$\Vc$ by (a) and (b).
\end{sm}

\begin{rmk} If $D$ is a finite group, $V$ and $H$ are subgroups of $D$
with double coset decomposition $D=\coprod_{P\in\clase} V\zeta_P
H$, then  $[D:H]=\sum_{P\in \clase} [V:V\bigcap \zeta_P
H\zeta^{-1}_P]$; if in addition $V\bigcap g Hg^{-1}=\{1\}$ for all
$g\in D$ then $\vert D\vert=\# \clase \vert V\vert\, \vert
H\vert$. See \cite[p. 76]{AM}.
\end{rmk}

We shall next analyze explicit examples of exact factors $\Hc$ of
$\Vc$ according to the equivalence relation $\approx_H$.

\medbreak
\subsection{Case 1. The equivalence relation $\approx_H$ is
 totally disconnected}

\

Here $\approx_H$ is the identity relation: if $P\approx_H Q$ then
$P=Q$ for any $P, Q\in\Pc$. Wide subgroupoids $\Hc$ with this
equivalence relation correspond to families $(H_P)_{P\in\Pc}$ of
subgroups  of $D$.

\begin{sm}\label{case1} Let $\Hc$ be as above.
Then the following statements are equivalent:
\begin{itemize}
\item[(i)] $\D = \Vc\Hc$ is an exact factorization.

\medbreak

\item[(ii)] $D=V H_P$ is an exact factorization for
any $P\in \Pc$.
\end{itemize}

\medbreak
To construct an explicit example of an exact factorization $\D =
\Vc\Hc$ with $\Vc$ connected and $\approx_H = \id$, we need a
finite group $D$, a subgroup $V$ and a family $(H_P)_{P\in\Pc}$ of
subgroups  of $D$ (thus the index set $\Pc$ is the basis of the
groupoid)  such that $H_P$ is an exact factor of $V$ for any
$P\in \Pc$.\end{sm}

Given a finite group $D$ and a subgroup $V$, we are thus faced to
the problem of finding all exact factors $H$ of $V$. We observe
that:

\medbreak (0).  Any conjugate of an exact factor of $V$ is again an exact factor of $V$.

\medbreak (1).  There exist a finite group $D$, a subgroup $V$ and
exact factors $H$ and $H^{'}$ such that $H\ncong H^{'}$. For
instance, $D=\Sy_4$, $V=\Sy_3$ (the subgroup that fixes 4),
$H=\langle(1234)\rangle\simeq \Z/(4)$,
$H^{'}=\langle(24)(13),(34)(12)\rangle\simeq \Z/(2)\oplus \Z/(2)$.

\medbreak (2). There exist a finite group $D$, a subgroup $V$ and
exact factors $H$ and $H^{'}$ with $H\simeq H^{'}$ but $H$ not
conjugate to $H^{'}$. For instance, $D=\Sy_n$, $n\ge 6$, $V=A_n$,
$H=\langle(12)\rangle$, $H^{'}=\langle(12)(34)(56)\rangle$.

\medbreak (3).  (Schur-Zassenhaus Theorem). If $D$ is a finite
group and $V\vartriangleleft D$ is a normal subgroup such that
$(|V|,[D:V])=1$, then $V$ admits exact factors, which are all
conjugate. The known proof of the conjugacy of the exact factors
relies on the Feit-Thompson Theorem (any group of odd order is
solvable).

\medbreak
(4). The list of all exact factorizations of $\Sy_n$ and $A_n$ is given in \cite{WW}.

\medbreak Let $\D = \Vc\Hc$ be an exact factorization with $\Vc$
connected and $\approx_H = \id$. Then, by Proposition \ref{coho2}
(i), the Kac exact sequence \eqref{kac} has the form

\begin{equation}\label{k1}
\begin{aligned}
0 &\longrightarrow  H^1(D, \ku^{\times})
\stackrel{\text{ res}}{\longrightarrow} \oplus_{P\in \Pc} H^1(H_P,
\ku^{\times}) \oplus H^1(V,
\ku^{\times})
\longrightarrow \Aut(\ku\,\Tc) \\
&\longrightarrow  H^2(D, \ku^{\times})
\stackrel{\text{ res}}{\longrightarrow}
\oplus_{P\in \Pc} H^2(H_P,\ku^{\times}) \oplus
H^2(V, \ku^{\times}) \longrightarrow  \Opext(\Vc,\Hc) \\
&\longrightarrow H^3(D, \ku^{\times})\stackrel{\text{
res }}{\longrightarrow} \oplus_{P\in \Pc} H^3(H_P,\ku^{\times}) \oplus
H^3(V, \ku^{\times})\longrightarrow \dots
\end{aligned}
\end{equation}

\begin{exa}
Let $m\in \Na$, $m\geq 5$. Let $D:=\Sy_{m}$, $V:=C_m = \langle
(1\dots m) \rangle$, $\Pc:=\{1 \dots m\}$, $H_i:=
\Sy_{m}^{(i)}:=\{\s\in\Sy_m : \s(i)=i\}$, $1\leq i\leq m$. The
groups $\Sy_{m}^{(i)}$ are conjugate to each other, indeed if
$\tau_{ij} = (ij)$ then
$\tau_{ij}\Sy_{m}^{(i)}\tau^{-1}_{ij}=\Sy_{m}^{(j)}$. We have
exact factorizations $\Sy_{m}=C_m\Sy_{m}^{(i)}$ for any $1\leq
i\leq m$, hence we have an exact factorization of groupoids
$\D=\Vc\Hc$.

\medbreak Assume now that $\ku=\C$. We rely on calculations done
in \cite{M2}. We claim that $\Opext(\Vc,\Hc)= \Z_2^{m-1}$. It is
known that $H^2(C_m,\ku^{\times})=0$ that
$H^2(\Sy_{m},\ku^{\times})=\Z_2$; that the restriction map $\text{
res }:H^2(\Sy_{m} ,\ku^{\times})\to H^2(\Sy_{m-1} ,\ku^{\times})$
is bijective \cite[p. 579]{M2}; and that $\text{ res }:H^3(\Sy_{m}
,\ku^{\times})\to H^3(\Sy_{m-1} ,\ku^{\times})\oplus H^3(C_m
,\ku^{\times})$ is an injective map. Hence the Kac exact sequence
\eqref{k1} gives $0 \longrightarrow  \Z_2 \longrightarrow  \Z_2^m
\longrightarrow  \Opext(\Vc,\Hc) \longrightarrow  0$ and the claim
follows.
\end{exa}

\subsection{Case 2. The equivalence relation $\approx_H$ is
connected}

\

Here both $\approx_V$ and $\approx_H$ are connected. A wide
subgroupoid $\Hc$ in this case corresponds to a subgroup $H$ of
$D$ and elements $(\lambda_P)_{P\in\Pc}\in H\backslash D$, such
that $\lambda_O \in H$.

\begin{sm}\label{caso2} Let $\Hc$ be a groupoid as above.
Then the following statements are equivalent.
\begin{itemize}

\item[(i)] $\D=\Vc\Hc$ is an exact factorization.

\item[(ii)] $D=\coprod_{r\in\Pc} V\lambda^{-1}_r H$
(thus $\# (V\backslash D/H)=\#\Pc$)
and  $V\bigcap g Hg^{-1}=\{1\}$ for any $g\in D$.
\end{itemize}

\medbreak To construct an explicit example of an exact
factorization $\D = \Vc\Hc$ with $\Vc$ and $\Hc$ connected, we
need a finite group $D$, two subgroups $V$ and $H$ such that
$V\cap g Hg^{-1}=\{1\}$ for any $g\in D$, and a section
$\zeta:V\backslash D/H \to D$ of the canonical projection $D \to
V\backslash D/H$ such that $\zeta(VH) \in H$. The base of the
groupoids $\D$, $\Vc$ and $\Hc$ is $\Pc:=V\backslash D/H$, and
$\zeta_P=\lambda^{-1}_P$ for any $P\in \Pc$. \end{sm}

Given a finite group $D$ and a subgroup $V$, we have to find
subgroups $H$ of $D$ such that $V$ intersects trivially all
conjugates of $H$. We observe that:

\medbreak
(1). If the orders $\vert V\vert$, $ \vert H\vert$ are relatively prime
then this condition is automatically fulfilled.

\medbreak (2).  If $V$ admits an exact factor $K$ and $H$ is any
subgroup of $K$ then $V$ intersects trivially all conjugates of
$H$.

\medbreak Let $\D = \Vc\Hc$ be an exact factorization with $\Vc$
and $\Hc$ connected. Then, by Proposition \ref{coho2} (i), the Kac
exact sequence \eqref{kac} has the form

\begin{equation}\label{kac2}
\begin{aligned}
0 &\longrightarrow  H^1(D, \ku^{\times})
\stackrel{\text{ res}}{\longrightarrow}  H^1(H,
\ku^{\times}) \oplus H^1(V,
\ku^{\times})
\longrightarrow \Aut(\ku\,\Tc) \\
&\longrightarrow  H^2(D, \ku^{\times})
\stackrel{\text{ res}}{\longrightarrow} H^2(H,
\ku^{\times}) \oplus
H^2(V, \ku^{\times}) \longrightarrow  \Opext(\Vc,\Hc) \\
&\longrightarrow H^3(D, \ku^{\times})\stackrel{\text{
res }}{\longrightarrow}  H^3(H,\ku^{\times}) \oplus
H^3(V, \ku^{\times})\longrightarrow \dots
\end{aligned}
\end{equation}

\begin{exa}Let $m, k, r\in\Na$ be such that $k,r\leq m$ and $k
> m-r$. Let $X\subseteq \{1,...,m\}$ be a subset of cardinal $r$.
Let $D:=\Sy_{m}$, $V:=C_k = \langle (12 \dots k) \rangle$, $H:=
\Sy^{X}_m:=\{\sigma\in\Sy_m: \sigma(x)=x \text{ for all } x\in X\}
\simeq \Sy_{m - r}$. Since $\sigma \Sy^{X}_m \sigma^{-1}=
\Sy^{\sigma(X)}_m$ for any $\sigma\in \Sy_m$, then $C_k\bigcap
\sigma \Sy^{X}_m \sigma^{-1} = \{e\}$ for any $\sigma\in \Sy_m$.
In this example $\# \Pc= \frac{n(n-1)...(n-r+1)}{k}$.
\end{exa}

\medbreak Let us assume that we are in conditions of Summary
\ref{caso2}. The set $VH\subseteq D$ is a subgroup of $D$ if and
only if $VH=HV$. If this is the case, we have an exact
factorization that allows us to consider the group $\Opext(V,H)$.
The following example shows that the isomorphism
$\Opext(V,H)\simeq \Opext(\Vc,\Hc)$ does not necessarily hold.

\medbreak

\begin{exa} In this example $\ku=\C$. Let $p,q,n\in\Na$, where $p$
and $q$ are relatively prime, and set $m=pqn$. Let $D:=\Z/(m)$,
$V:=\Z/(p)$, $H:=\Z/(q)$. Let $\D, \Vc$ and $\Hc$ the
corresponding connected groupoids associated to $D,V$    and $H$
respectively, with $\Pc=V\backslash D/H$ of cardinal $n$.

\medbreak

We claim that $\Opext(\Vc,\Hc)\simeq \Z/(n)$. Indeed, since
$H^2(\Z/(p), \C^{\times})=H^2(\Z/(q), \C^{\times})=0$ then
$\Opext(\Vc,\Hc)\simeq \Ker(\text{res}:H^3(\Z/(m),\C^{\times})\to
H^3(\Z/(p),\C^{\times})\oplus H^3(\Z/(q),\C^{\times}))$ by
\eqref{kac2}. Now, it is known that $ H^3(\Z/(r),
\C^{\times})\simeq \Z/(r)$ for any $r\in\Na$ \cite[p. 61]{AM}. And
it is not difficult to see that the restriction map $\text{res}:
\Z/(m)\to \Z/(p)\oplus \Z/(q)$ via these isomorphisms is the
canonical projection, and thus  $\Opext(\Vc,\Hc)\simeq \Z/(n)\not
\simeq \Opext(V,H)$ unless $n=1$.
\end{exa}

\subsection{Case 3. Equivalence relations $\approx_H$ with two classes}

\

We assume here that $\approx_H$ has two equivalence classes:
$\{O\}$, and $\Pc-\{O\}$. We fix $\widetilde{O}\in\Pc-\{O\}$. A
wide subgroupoid $\Hc$ in this case corresponds to a pair of
subgroups $H_1$, $H_2$ and a family $(\lambda_P)_{P\in\Pc-\{O\}}$
in $D$ such that $\lambda_{\widetilde{O}}\in H_2$.

\begin{sm}\label{case3} Let $\Hc$ be a groupoid as above. The
following statements are equivalent.
\begin{itemize}
\item[(i)] $\D=\Vc\Hc$ is an exact factorization,
\item[(ii)] \begin{itemize}\item[(a)] $D=VH_1$ is an exact
factorization,
\item[(b)] $D=\coprod_{P\in\Pc-\{O\}} V\lambda^{-1}_P H_2$
and $V\bigcap g H_2 g^{-1}=\{1\}$ for any $g\in D$. \end{itemize}
\end{itemize}

By Proposition \ref{coho2} (i) the Kac exact sequence has the form
\begin{equation}\label{kac3}
\begin{aligned}
0 &\longrightarrow  H^1(D, \ku^{\times}) \stackrel{\text{
res}}{\longrightarrow} H^1(H_1, \ku^{\times}) \oplus H^1(H_2,
\ku^{\times}) \oplus  H^1(V, \ku^{\times})
\longrightarrow \Aut(\ku\,\Tc) \\
&\longrightarrow  H^2(D, \ku^{\times}) \stackrel{\text{
res}}{\longrightarrow} H^2(H_1, \ku^{\times}) \oplus H^2(H_2,
\ku^{\times})  \oplus H^2(V, \ku^{\times})
\longrightarrow  \Opext(\Vc,\Hc) \\
&\longrightarrow H^3(D, \ku^{\times})\stackrel{\text{ res
}}{\longrightarrow} H^3(H_1,\ku^{\times}) \oplus
H^3(H_2,\ku^{\times}) \oplus H^3(V, \ku^{\times})\longrightarrow
\dots
\end{aligned}
\end{equation}
\end{sm}

Examples in this case are obtained combining the examples in the
previous cases. We note a general way of obtaining collections as
in Summary \ref{case3}. If  $D=VH_1$ is an exact factorization and
$H_2$ is a subgroup of $H_1$, then $V\bigcap g H_2 g^{-1}=\{1\}$
for any $g\in D$. Let us discuss an explicit example.

\begin{exa} Let $H_2$ be any group of order $n$, considered as a subgroup
of $\Sy_n$ via, say, the left regular action on itself. Set
$D=\Sy_{n+1}$, $V=<(12\dots n+1)>\simeq C_{n+1}$, $H_1=\Sy_n$, an
exact factor of  $V$. By \eqref{kac3}, we have
\[\Opext(\Vc,\Hc) \cong H^2(H_2,
\ku^{\times}) / \Im\left(\text{res}:H^2(\Sy_{n+1},
\ku^{\times})\to H^2(H_2, \ku^{\times})\right).\]
\end{exa}

There are also examples which are not of this form.

\begin{exa} Let $n\in\Na$ such that $n=rs$ with $r,s\neq 1$, $r,s\in \Na$.
Let $D=\Sy_n$, $V= \Sy_{n-1}$, $H_1=<(1,...,n)>$ and
$H_2=<\sigma>$, where $\sigma = (1 \dots r)(r+1 \dots 2r) \dots
(rs - r + 1 \dots rs)$. Then $H_1$ is an exact factor of $V$, $V$
intersects trivially any conjugate of $H_2$ and $\sigma \notin
H_1$. In this case $\# \Pc= s+1$.
\end{exa}

\bibliographystyle{amsalpha}

\end{document}